\documentclass{amsart}
\usepackage[all]{xy}
\usepackage{amscd,amsthm,amssymb,amsfonts,amsmath,euscript}

\theoremstyle{plain}
\newtheorem{thm}{Theorem}[section]
\newtheorem{lm}[thm]{Lemma}
\newtheorem{prop}[thm]{Proposition}
\newtheorem{cor}[thm]{Corollary}

\theoremstyle{definition}
\newtheorem{defn}[thm]{Definition}

\theoremstyle{remark}
\newtheorem{remark}[thm]{Remark}

\newcommand{\nc}{\newcommand}

\numberwithin{equation}{section}

\def\makeop#1{\expandafter\def\csname#1\endcsname
  {\mathop{\rm #1}\nolimits}\ignorespaces}
\makeop{Hom}   \makeop{End}   \makeop{Aut}   \makeop{Isom}  \makeop{Pic}
\makeop{Gal}   \makeop{ord}   \makeop{Char}  \makeop{Div}   \makeop{Lie}
\makeop{PGL}   \makeop{Corr}  \makeop{PSL}   \makeop{sgn}   \makeop{Spf}
\makeop{Spec}  \makeop{Tr}    \makeop{Nr}    \makeop{Fr}    \makeop{disc}
\makeop{Proj}  \makeop{supp}  \makeop{ker}   \makeop{im}    \makeop{dom}
\makeop{coker} \makeop{Stab}  \makeop{SO}    \makeop{SL}    \makeop{SL}
\makeop{Cl}    \makeop{cond}  \makeop{Br}    \makeop{inv}   \makeop{rank}
\makeop{id}    \makeop{Fil}   \makeop{Frac}  \makeop{GL}    \makeop{SU}
\makeop{Trd}   \makeop{Sp}    \makeop{Tr}    \makeop{Trd}   \makeop{diag}
\makeop{Res}   \makeop{ind}   \makeop{depth} \makeop{Tr}    \makeop{st}
\makeop{Ad}    \makeop{Int}   \makeop{tr}    \makeop{Sym}   \makeop{can}
\makeop{length}\makeop{SO}    \makeop{torsion} \makeop{GSp} \makeop{pr}
\def\makebb#1{\expandafter\def
  \csname bb#1\endcsname{{\mathbb{#1}}}\ignorespaces}
\def\makebf#1{\expandafter\def\csname bf#1\endcsname{{\bf
      #1}}\ignorespaces}
\def\makegr#1{\expandafter\def
  \csname gr#1\endcsname{{\mathfrak{#1}}}\ignorespaces}
\def\makescr#1{\expandafter\def
  \csname scr#1\endcsname{{\EuScript{#1}}}\ignorespaces}
\def\makecal#1{\expandafter\def\csname cal#1\endcsname{{\mathcal
      #1}}\ignorespaces}

\def\doLetters#1{#1A #1B #1C #1D #1E #1F #1G #1H #1I #1J #1K #1L #1M
                 #1N #1O #1P #1Q #1R #1S #1T #1U #1V #1W #1X #1Y #1Z}
\def\doletters#1{#1a #1b #1c #1d #1e #1f #1g #1h #1i #1j #1k #1l #1m
                 #1n #1o #1p #1q #1r #1s #1t #1u #1v #1w #1x #1y #1z}
\doLetters\makebb   \doLetters\makecal  \doLetters\makebf
\doLetters\makescr
\doletters\makebf   \doLetters\makegr   \doletters\makegr
     \def\qed{\qedmark\medbreak}%
\def\qedmark{{\enspace\vrule height 6pt width 5pt depth 1.5pt}}%
    
   \def\Ga{{{\bbG}_{\rm a}}}

\long\def\MSC#1\EndMSC{\def\arg{#1}\ifx\arg\empty\relax\else
     {\par\narrower\noindent%
     2000 Mathematics Subject Classification: #1\par}\fi}

\long\def\KEY#1\EndKEY{\def\arg{#1}\ifx\arg\empty\relax\else
        {\par\narrower\noindent Keywords and Phrases: #1\par}\fi\TSkip}

\normalsize

\makeop{Bl}

\def\Spec{{\rm Spec}}

\def\Fpbar{\overline{\bbF}_p}
\def\Fp{{\bbF}_p}
\def\Fq{{\bbF}_q}

\def\Qbar{\overline{\bbQ}}

\newcommand{\Z}{\mathbb Z}
\newcommand{\Q}{\mathbb Q}

\newcommand{\C}{\mathbb C}
\newcommand{\F}{\mathbb F}

\newcommand{\npr}{\noindent }

\newcommand{\<}{\langle}   
\renewcommand{\>}{\rangle} 

\nc{\embed}{\hookrightarrow}

\newcommand{\ch}{characteristic }
\newcommand{\ac}{algebraically closed }
\newcommand{\dieu}{Dieudonn\'{e} }

\nc{\ol}{\overline}
\nc{\wt}{\widetilde}
\nc{\opp}{\mathrm{opp}}
\nc{\ul}{\underline}

\makeop{Ram}
\makeop{Rep}

\begin{document}
\renewcommand{\thefootnote}{\fnsymbol{footnote}}
\setcounter{footnote}{-1}
\title{The supersingular loci and mass formulas on Siegel modular varieties}
\author{Chia-Fu Yu}
\address{
Institute of Mathematics \\
Academia Sinica  \\
128 Academia Rd.~Sec.~2, Nankang\\
Taipei, Taiwan\\
and NCTS (Taipei Office)}
\email{chiafu@math.sinica.edu.tw}

\date{December 22, 2006. The research is partially supported by NSC
  95-2115-M-001-004. 
}

\begin{abstract}
We describe the supersingular locus of the Siegel 3-fold with a 
parahoric level structure. We also study its higher dimensional
generalization. Using this correspondence and a deep result of Li and
Oort, we evaluate the number of irreducible components
of the supersingular locus of the Siegel moduli space $\calA_{g,1,N}$ 
for arbitrary $g$.
\end{abstract}
\maketitle


\section{Introduction}
\label{sec:01}

In this paper we discuss some extensions of works of Katsura and
Oort \cite{katsura-oort:surface}, and of Li and Oort \cite{li-oort} on
the supersingular locus of a mod $p$ Siegel modular variety.   

Let $p$ be a rational prime number, $N\ge 3$ a prime-to-$p$
positive integer. We choose a primitive $N$-th root of unity $\zeta_N$ in
$\Qbar\subset \C$ and an embedding $\Qbar\embed \Qbar_p$.
Let $\calA_{g,1,N}$ be the moduli space over $\Z_{(p)}[\zeta_N]$ of
$g$-dimensional principally polarized abelian varieties
$(A,\lambda,\eta)$ with a symplectic level-$N$ structure (See
Subsection~\ref{sec:21}).

Let $\calA_{2,1,N,(p)}$ be the cover of $\calA_{2,1,N}$ which
parametrizes isomorphism classes of objects $(A,\lambda,\eta,
H)$, where 
$(A,\lambda,\eta)$ is an object in $\calA_{2,1,N}$ and $H\subset A[p]$
is a finite flat subgroup scheme of rank $p$.
It is known that the moduli scheme $\calA_{2,1,N,(p)}$ has semi-stable
reduction and
the reduction $\calA_{2,1,N,(p)}\otimes \Fpbar$ modulo $p$ has two
irreducible components.
Let $\calS_{2,1,N,(p)}$
(resp. $\calS_{2,1,N}$) denote the
supersingular locus of the moduli space $\calA_{2,1,N,(p)}\otimes
\Fpbar$ (resp. $\calA_{2,1,N}\otimes \Fpbar$). Recall that an abelian
variety $A$ in \ch $p$ 
  is called {\it supersingular} if it is isogenous to a product of
  supersingular elliptic curves over an \ac field $k$; it is called
  {\it superspecial} if it is isomorphic to a product of supersingular
  elliptic curves over $k$.

The supersingular locus
$\calS_{2,1,N}$ of the Siegel 3-fold is studied
in Katsura and Oort \cite{katsura-oort:surface}. We summarize the main
results for $\calS_{2,1,N}$ (the local results
obtained earlier in Koblitz \cite{koblitz:thesis}):

\begin{thm} \label{11} \ 

{\rm (1)} The scheme $\calS_{2,1,N}$ is equi-dimensional and each
  irreducible component is isomorphic to $\bfP^1$. 

{\rm (2)} The scheme $\calS_{2,1,N}$  has
  \begin{equation}
    \label{eq:11}
  |\Sp_4(\Z/N\Z)|\frac{(p^2-1)}{5760}     
  \end{equation}
  irreducible components.

{\rm (3)} The singular points of $\calS_{2,1,N}$ are exactly the
    superspecial points and there are 
    \begin{equation}
      \label{eq:12}
      |\Sp_4(\Z/N\Z)|\frac{(p-1)(p^2+1)}{5760}
    \end{equation}
    of them. Moreover, at each singular point there are $p+1$ irreducible
    components passing through and intersecting transversely.
\end{thm}
\begin{proof}
  See Koblitz \cite[p.193]{koblitz:thesis} and Katsura-Oort
  \cite[Section 2, Theorem 5.1, Theorem 5.3]{katsura-oort:surface}. \\
\end{proof}

In this paper we extend their results to $\calS_{2,1,N,(p)}$. We show

\begin{thm} \label{12} \

{\rm (1)} The scheme $\calS_{2,1,N,(p)}$ is equi-dimensional and
  each irreducible component is isomorphic to
  $\bfP^1$.

{\rm (2)} The scheme $\calS_{2,1,N,(p)}$  has
  \begin{equation}\label{eq:13}
 |\Sp_4(\Z/N\Z)|\cdot \frac{(-1)\zeta(-1)\zeta(-3)}{4} \left
  [(p^2-1)+(p-1)(p^2+1)\right ]  
  \end{equation}
  irreducible components, where $\zeta(s)$ is the Riemann zeta
  function.

{\rm (3)} The scheme $\calS_{2,1,N,(p)}$ has only ordinary double
singular points and there are
\begin{equation}
  \label{eq:14}
 |\Sp_4(\Z/N\Z)|\cdot \frac{(-1)\zeta(-1)\zeta(-3)}{4}
(p-1)(p^2+1)(p+1)  
\end{equation}
of them. 

{\rm (4)} The natural morphism $\calS_{2,1,N,(p)}\to \calS_{2,1,N}$ contracts
\begin{equation}
  \label{eq:15}
  |\Sp_4(\Z/N\Z)|\cdot \frac{(-1)\zeta(-1)\zeta(-3)}{4}(p-1)(p^2+1)
\end{equation}
projective lines onto the superspecial points of
$\calS_{2,1,N}$. 
\end{thm}
\begin{remark}(1)
By the basic fact that
\[ \zeta(-1)=\frac{-1}{12}\quad \text{and}\quad
\zeta(-3)=\frac{1}{120}, \]
the number (\ref{eq:15}) (of the vertical components)
equals the number (\ref{eq:12}) (of superspecial points), and the
number (\ref{eq:13}) (sum of vertical and horizontal components)
equals the sum of the numbers (\ref{eq:11}) (of irreducible components) and
(\ref{eq:12}) (of superspecial points). Thus, the set of 
horizontal irreducible components of $\calS_{2,1,N,(p)}$ is in
bijection with the set of irreducible components of $\calS_{2,1,N}$

(2) Theorem~\ref{12} (4) says that $\calS_{2,1,N,(p)}$ is a
    ``desingularization'' or a ``blow-up'' of $\calS_{2,1,N}$ at
    the singular points. Strictly speaking, the desingularization of
    $\calS_{2,1,N}$ is its normalization, which is the (disjoint)
    union of horizontal components of
    $\calS_{2,1,N,(p)}$. The vertical components of $\calS_{2,1,N,(p)}$
    should be the exceptional divisors of the blowing up of a suitable
    ambient surface of $\calS_{2,1,N}$ at superspecial points.  
\end{remark}

In the proof of Theorem~\ref{12} (Section~\ref{sec:04}) we see that 
\begin{itemize}
\item the set of certain superspecial points (the set $\Lambda$ in
  Subsection~\ref{sec:41}) in $\calS_{2,p,N}$ (classifying degree $p^2$
  polarized supersingular abelian surfaces) is in bijection with the
  set of irreducible components of $\calS_{2,1,N}$, and
\item the set of superspecial points in
$\calS_{2,1,N}$ is in bijection with the set of irreducible
components of $\calS_{2,p,N}$, furthermore
\item the supersingular locus $\calS_{2,1,N,(p)}$ provides  the explicit
  link of this duality as a correspondence that performs simply through
  the ``blowing-ups'' and ``blowing-downs''.
\end{itemize}
In the second part of this paper we attempt to generalize a similar
picture to higher (even) dimensions. \\

Let $g=2D$ be an even positive integer. 
Let $\calH$ be the moduli space over $\Z_{(p)}[\zeta_N]$ which 
parametrizes equivalence classes of objects $(\varphi:\ul A_1\to
\ul A_2)$,  where 
\begin{itemize}
\item $\ul A_1=(A_1,\lambda_1,\eta_1)$ is an object in
$\calA_{g,1,N}$,
\item $\ul A_2=(A_2,\lambda_2,\eta_2)$ is an object in
$\calA_{g,p^D,N}$, and
\item  $\varphi:A_1\to A_2$ is an
isogeny of degree $p^D$
satisfying $\varphi^*\lambda_2=p\lambda_1$ and
$\varphi_*\eta_1=\eta_2$.
\end{itemize}
 
The moduli space $\calH$ with natural projections gives the
following correspondence:
\[ \xymatrix{
 &  \calH \ar[ld]_{\pr_1} \ar[rd]^{\pr_2} & \\
\calA_{g,1,N} & & \calA_{g,p^D,N}.
} \]

Let $\calS$ be the supersingular locus of $\calH\otimes \Fpbar$, which
is the reduced closed subscheme consisting of supersingular points
(either $A_1$ or $A_2$ is supersingular, or equivalently both are so).


In the special case where $g=2$, $\calH$ is isomorphic to
$\calA_{2,1,N,(p)}$, and $\calS \simeq \calS_{2,1,N,(p)}$ under this
isomorphism (See Subsection~\ref{s:45}). 

As the second main result of this paper, we obtain

\begin{thm}\label{13}
  Let $C$ be the number of irreducible
  components of $\calS_{g,1,N}$. Then
\[ C=|\Sp_{2g}(\Z/N\Z)|\cdot
  \frac{(-1)^{g(g+1)/2}}{2^g} \left \{ \prod_{i=1}^g \zeta(1-2i)
  \right \}\cdot L_p\, , \]
where 
\[ L_p=
\begin{cases}
  \prod_{i=1}^{g}\left\{(p^i+(-1)^i\right \} & \text{if $g$ is odd;}
  \\
  \prod_{i=1}^{D}(p^{4i-2}-1) & \text{if $g=2D$ is even.} 
\end{cases}\]
\end{thm}

In the special case where $g=2$, Theorem~\ref{13} recovers
 Theorem~\ref{11} (2). 

We give the idea of the proof. 
Let $\Lambda_{g,1,N}$ denote the set of superspecial
(geometric) points in $\calA_{g,1,N}\otimes \Fpbar$.
For $g=2D$ is even, let $\Lambda^*_{g,p^{D},N}$
denote the set of superspecial (geometric) points $(A,\lambda,\eta)$
in $\calA_{g,p^D,N}\otimes \Fpbar$ satisfying $\ker \lambda=A[F]$,
where $F:A\to A^{(p)}$ is the relative Frobenius morphism on $A$. 
These are finite sets and each member is defined over $\Fpbar$.
By a result of Li and Oort~\cite{li-oort} (also see
Section~\ref{sec:05}), we know 
\[ C=
\begin{cases}
 |\Lambda_{g,1,N}| &  \text{if $g$ is odd;} \\
 |\Lambda^*_{g,p^{D},N}| &  \text{if $g$ is even.}
\end{cases} \]
One can use the geometric mass formula due to 
Ekedahl~\cite{ekedahl:ss} and some others (see Section~\ref{sec:03})  
to compute $|\Lambda_{g,1,N}|$. Therefore, it 
remains to compute $|\Lambda^*_{g,p^{D},N}|$ when $g$ is even. 
We restrict the
correspondence $\calS$ to the product
$\Lambda_{g,1,N}\times\Lambda^*_{g,p^{D},N}$ of superspecial
points, and compute certain special points in $\calS$. This gives us 
relation between $\Lambda^*_{g,p^{D},N}$ and $\Lambda_{g,1,N}$. 
See Section~\ref{sec:06} for details. 

Theorem~\ref{13} tells us how the number $C=C(g,N,p)$
varies when $p$ varies. For another application,
one can use this result to compute the 
dimension of the space of Siegel cusp forms of certain level at $p$ by
the expected Jacquet-Langlands correspondence for symplectic
groups. As far as the author 
knows, the latter for general $g$ is not available yet in the
literature.   

The paper is organized as follows. In Section~2, we recall the basic
definitions and properties of the Siegel moduli spaces and
supersingular abelian varieties. In Section~3, we state the mass
formula for superspecial principally polarized abelian varieties due
to Ekedahl (and some others). The proof of Theorem~\ref{12} is given
in Section 4. In Section 5, we describe the results of Li and Oort on
irreducible components of the supersingular locus. In Section 6,
we introduce a correspondence and use this to evaluate the number of
irreducible components of the supersingular locus. \\

\npr {\it Acknowledgments.} The author is grateful to Katsura, Li and
Oort for their inspiring papers on which the present work relies heavily. 
He also thanks C.-L.~Chai for his encouragements on this subject for
years, and thanks T.~Ibukiyama for his interest of the present work.  
He is indebted to the referee for careful reading and helpful comments
that improve the presentation of this paper.




\section{Notation and preliminaries}
\label{sec:02}

\subsection{}\label{sec:21}
Throughout this paper we fix a rational prime $p$ and a prime-to-$p$
positive integer $N\ge 3$. Let $d$ be a positive integer with $(d,N)=1$.
We choose a primitive $N$-th root of unity $\zeta_N$ in
$\Qbar\subset \C$ and an embedding $\Qbar\embed \Qbar_p$. The element
$\zeta_N$ gives rise to a trivialization
$\Z/N\Z\simeq \mu_N$ over any $\Z_{(p)}[\zeta_N]$-scheme.
For a polarized abelian variety $(A,\lambda)$ of degree $d^2$, a
full symplectic level-$N$
structure with respect to the choice $\zeta_N$ is an isomorphism
\[ \eta: (\Z/N\Z)^{2g}\simeq A[N] \]
such that the following diagram commutes
\[
\begin{CD}
   (\Z/N\Z)^{2g}\times  (\Z/N\Z)^{2g} @>(\eta, \eta)>> A[N]\times A[N]
   \\
   @V\< \, , \,\,\> VV @V e_\lambda VV \\
   \Z/N\Z @>\zeta_N >> \mu_N, \\
\end{CD}
\]
where $\<\, ,\>$ is the standard non-degenerate alternating form on 
$(\Z/N\Z)^{2g}$ and $e_\lambda$ is the Weil pairing induced by the
polarization $\lambda$. 

Let $\calA_{g,d,N}$ denote the moduli space over $\Z_{(p)}[\zeta_N]$ of
$g$-dimensional polarized abelian varieties $(A,\lambda,\eta)$ of
degree $d^2$ with a
full symplectic level $N$ structure with respect to $\zeta_N$. 
Let $\calS_{g,d,N}$  denote the supersingular locus of the reduction
$\calA_{g,d,N}\otimes \Fpbar$ modulo $p$, which is the
closed reduced subscheme of $\calA_{g,d,N}\otimes \Fpbar$ consisting of
supersingular points in $\calA_{g,d,N}\otimes \Fpbar$. 
Let $\Lambda_{g,d,N}$ denote the set of superspecial
(geometric) points in $\calS_{g,d,N}$;
this is a finite set and every member is defined over $\Fpbar$. \\

For a scheme $X$ of finite type over a field $K$, we denote by
$\Pi_0(X)$ the set of geometrically irreducible components of $X$. 

Let $k$ be an \ac field of \ch $p$.

\subsection{} Over a ground field $K$ of \ch $p$, denote by $\alpha_p$
the finite group scheme of rank $p$ that is the kernel of the
Frobenius endomorphism from the additive group $\Ga$ to itself. One has
\[ \alpha_p=\Spec K[X]/X^p, \quad m(X)=X\otimes 1+1\otimes X, \]
where $m$ is the group law.

By definition, an elliptic curve $E$ over $K$ is called {\it
  supersingular} if $E[p](\ol K)=0$. An abelian variety $A$ over $K$
  is called {\it supersingular} if it is isogenous to a product of
  supersingular elliptic curves over $\ol K$; $A$ is called 
  {\it superspecial} if it is isomorphic to a product of supersingular
  elliptic curves over $\ol K$.

For any abelian variety $A$ over $K$ where $K$ is perfect, the {\it
$a$-number} of $A$ is defined by  
\[ a(A):=\dim_K \Hom(\alpha_p,A). \]

The following interesting results are well-known; see Subsection 1.6
of \cite{li-oort} for a detail discussion. 
\begin{thm}[Oort]
  If $a(A)=g$, then $A$ is superspecial.
\end{thm}

\begin{thm}[Deligne, Ogus, Shioda]
  For $g\ge 2$, there is only one
  $g$-dimensional superspecial abelian variety up to
  isomorphism over $k$.
\end{thm}

\section{The mass formula}
\label{sec:03}
Let $\Lambda_g$ denote the set of isomorphism classes of
$g$-dimensional principally polarized superspecial abelian varieties
over $\Fpbar$. Write
\[ M_g:=\sum_{(A,\lambda)\in\Lambda_g}
  \frac{1}{|\Aut(A,\lambda)|} \]
for the mass attached to $\Lambda_g$. The following mass formula is
due to Ekedahl \cite[p.159]{ekedahl:ss} and Hashimoto-Ibukiyama
\cite[Proposition 9]{hashimoto-ibukiyama:classnumber}.


\begin{thm}\label{31} Notation as above. One has
    \begin{equation}
  M_g=
  \frac{(-1)^{g(g+1)/2}}{2^g} \left \{ \prod_{k=1}^g \zeta(1-2k)
  \right \}\cdot \prod_{k=1}^{g}\left\{(p^k+(-1)^k\right \}.
  \end{equation}
\end{thm}
Similarly, we set
\[ M_{g,1.N}:=\sum_{(A,\lambda,\eta)\in\Lambda_{g,1,N}}
  \frac{1}{|\Aut(A,\lambda,\eta)|}. \]
\begin{lm}\label{32} We have
   $M_{g,1,N}=|\Lambda_{g,1,N}|=|\Sp_{2g}(\Z/N\Z)|\cdot M_g.$
\end{lm}
\begin{proof}
  The first equality follows from a basic fact that $(A,\lambda,\eta)$
  has no non-trivial automorphism. The proof of the second equality is
  elementary; see Subsection 4.6 of \cite{yu:mass_hb}. 
\end{proof}
\begin{cor}\label{33} One has
 \[ |\Lambda_{2,1,N}|=|\Sp_4(\Z/N\Z)|\cdot
 \frac{(-1)\zeta(-1)\zeta(-3)}{4}(p-1)(p^2+1). \]
\end{cor}

\section{Proof of Theorem~\ref{12}}
\label{sec:04}
\subsection{}\label{sec:41} 
In this section we consider the case where $g=2$. Let 
\[ \Lambda:=\{(A,\lambda,\eta)\in \calS_{2,p,N};\ \ker \lambda\simeq
\alpha_p\times \alpha_p \}. \]
Note that every member $\ul A$ of $\Lambda$ is superspecial (because
$A\supset \alpha_p\times \alpha_p$), that is, 
$\Lambda \subset \Lambda_{2,p,N}$. 
For a point $\xi$ in $\Lambda$, consider the space $S_\xi$ which
parametrizes the isogenies of degree $p$
\[ \varphi: (A_\xi,\lambda_\xi,\eta_\xi)\to \ul A=(A,\lambda,\eta) \]
which preserve the polarizations and level structures. Let
\[ \psi_\xi: S_\xi\to \calS_{2,1,N} \]
be the morphism which sends $(\varphi:\xi\to \ul A)$ to $\ul A$. Let
$V_\xi\subset \calS_{2,1,N}$ be the image of $S_\xi$ under
$\psi_\xi$.

The following theorem is due to Katsura and Oort
\cite[Theorem 2.1 and Theorem 5.1]{katsura-oort:surface}:

\begin{thm}[Katsura-Oort]
  Notation as above. The map $\xi\mapsto V_\xi$ gives rise to a bijection
$\Lambda\simeq \Pi_0(\calS_{2,1,N})$ and one has
 \[ |\Lambda|=|\Sp_4(\Z/N\Z)|{(p^2-1)}/{5760}. \]
\end{thm}

We will give a different way of evaluating $|\Lambda|$ that is 
based on the geometric mass formula (see Corollary~\ref{44}).

\subsection{} Dually we can consider the space $S'_\xi$ for each
$\xi\in \Lambda$
that parametrizes the isogenies of degree $p$
\[ \varphi':\ul A=(A,\lambda,\eta)\to \xi=(A_\xi,\lambda_\xi,\eta_\xi), \]
with $\ul A\in \calA_{2,1,N}\otimes \Fpbar$, such that
$\varphi'_*\eta=\eta_\xi$ and 
$\varphi'^*\lambda_\xi=p\,\lambda$.  Let
\[ \psi'_\xi: S'_\xi\to \calS_{2,1,N} \]
be the morphism which sends $(\varphi':\ul A\to \xi)$ to $\ul A$. Let
$V'_\xi\subset \calS_{2,1,N}$ be the image of $S'_\xi$ under
$\psi'_\xi$.

For a degree $p$ isogeny $(\varphi: \ul A_1\to \ul A_2)$ with
  $\ul A_2$ in $\calA_{2,1,N}$, $\varphi^*\lambda_2=\lambda_1$ and
  $\varphi_*\eta_1=\eta_2$, we define  \[(\varphi: \ul A_1\to \ul
  A_2)^*=(\varphi': \ul A_2'\to \ul A_1'),\] where
  $\varphi'=\varphi^t$ and
\[ \ul  A_2'=(A_2^t,\lambda_2^{-1},\lambda_2\circ \eta_2), \]
\[ \ul  A_1'=(A_1^t, p\,\lambda_1^{-1},\lambda_1\circ \eta_1). \]
Note that $\varphi'_* \eta_2'=\eta_1'$ as we have the commutative
diagram:
\[
\begin{CD}
  (\Z/N\Z)^{4} @>\eta_2>> A_2[N] @>\lambda_2>> A_2^t[N] \\
    @VV=V  @AA\varphi A  @VV\varphi^t V     \\
   (\Z/N\Z)^{4} @>\eta_1>> A_1[N] @>\lambda_1>> A_1^t[N]. \\
\end{CD}\]
If $\ul A_1 \in \Lambda$, then $\ul A_1'$ is also in
$\Lambda$. Therefore,
the map $\xi\mapsto V'_\xi$ also gives rise to a bijection
$\Lambda\simeq \Pi_0(\calS_{2,1,N})$.

\subsection{} We use the classical contravariant \dieu theory. 
We refer the reader to
Demazure \cite{demazure} for a basic account of this theory. Let $K$
be a perfect field of \ch $p$, $W:=W(K)$ the ring of Witt vectors over
$K$, $B(K)$ the fraction field of $W(K)$. Let $\sigma$ be the
Frobenius
map on $B(K)$. A quasi-polarization on a \dieu module $M$ here is a
non-degenerate (meaning of non-zero discriminant) alternating pairing
\[ \<\, ,\>:M\times M\to B(K), \]
such that $\<Fx,y\>=\<x,Vy\>^\sigma$ for $x, y\in M$ and
$\<M^t,M^t\>\subset W$. Here we regard  
the dual $M^t$ of $M$ as a \dieu submodule in $M\otimes B(K)$ using
the pairing. A quasi-polarization is called {\it separable} if $M^t=M$.
Any polarized abelian variety $(A,\lambda)$ over $K$ naturally gives
rise to a quasi-polarized \dieu module. The induced 
quasi-polarization is separable if and only if $(p,\deg \lambda)=1$.

Recall (Subsection~\ref{sec:21}) that $k$ denotes an \ac field of \ch $p$.  

\begin{lm}\label{42} \ 

  {\rm (1)} Let $M$ be a separably quasi-polarized superspecial \dieu
  module over $k$ of rank $4$. Then there exists a basis
  $f_1,f_2,f_3,f_4$ for $M$ over $W:=W(k)$ 
  such that
  \[ F f_1=f_3, Ff_3=pf_1, \quad Ff_2=f_4, Ff_4=pf_2 \]
    and the non-zero pairings are
  \[ \<f_1,f_3\>=-\<f_3,f_1\>=\beta_1,\quad 
  \<f_2,f_4\>=-\<f_4,f_2\>=\beta_1,\] 
    where $\beta_1\in W(\F_{p^2})^\times$ with
  $\beta_1^\sigma=-\beta_1$.

  {\rm (2)} Let $\xi$ be a point in $\Lambda$, and let $M_{\xi}$ be
  the \dieu module of $\xi$. Then there is a $W$-basis 
  $e_1, e_2, e_3, e_4$ for $M_\xi$ such that
  \[ Fe_1=e_3,\quad Fe_2=e_4,\quad Fe_3=pe_1, \quad Fe_4=pe_2, \]
  and the non-zero pairings are
  \[  \<e_1,e_2\>=-\<e_2,e_1\>=\frac{1}{p}, \quad
  \<e_3,e_4\>=-\<e_4,e_3\>=1. \] 
\end{lm}
\begin{proof}
  (1) This is a special case of Proposition 6.1 of \cite{li-oort}.

  (2) By Proposition 6.1 of \cite{li-oort}, $(M_\xi,\<\,,\>)$
  either is indecomposable or decomposes into a product of two
  quasi-polarized supersingular \dieu modules of rank 2. In the
  indecomposable case, one can choose such a basis $e_i$ for
  $M_\xi$. Hence it remains to show that $(M_\xi,\<\,,\>)$ is
  indecomposable. 
  Let $(A_\xi[p^\infty],\lambda_\xi)$ be the associated polarized
  $p$-divisible group. Suppose it decomposes into
  $(H_1,\lambda_1)\times (H_2,\lambda_2)$. Then the kernel of
  $\lambda$ is isomorphic to 
  $E[p]$ for a supersingular elliptic curve $E$. Since $E[p]$ is a
  nontrivial extension of $\alpha_p$ by $\alpha_p$, one gets
  contradiction. This completes the proof. \qed
\end{proof}

\subsection{} \label{s:44} Let $(A_0,\lambda_0)$ be a superspecial
principally polarized abelian surface and $(M_0,\<\,,\>_0)$ be the
associated 
\dieu module. Let $\varphi':(A_0,\lambda_0)\to (A,\lambda)$ be an
isogeny of degree $p$ with $\varphi'^* \lambda=p\,\lambda_0$.
Write $(M,\<\, ,\>)$ for the \dieu module of $(A,\lambda)$. Choose a basis
$f_1,f_2,f_3,f_4$ for $M_0$ as in Lemma~\ref{42}. We have the inclusions
  \[ (F,V)M_0\subset M \subset M_0. \]
Modulo $(F,V)M_0$, a module $M$ corresponds a one-dimensional subspace
$\ol M$ in $\ol M_0:=M_0/(F,V)M_0$. As $\ol M_0=k<f_1,f_2>$, $\ol M$
is of the form
\[ \ol M=k<af_1+bf_2>, \quad [a:b]\in \bfP^1(k). \]

The following result is due to Moret-Bailly
\cite[p.138-9]{moret-bailly:p1}. We include a proof for the reader's
convenience.

\begin{lm}\label{43}
Notation as above, $\ker
\lambda\simeq \alpha_p\times \alpha_p$ if and only if the
corresponding point $[a:b]$ satisfies
$a^{p+1}+b^{p+1}=0$. Consequently, there are $p+1$ isogenies $\varphi'$
so that $\ker \lambda\simeq \alpha_p\times \alpha_p$.
\end{lm}
\begin{proof}
  As $\varphi'^*\lambda=p\lambda_0$, we have 
$\<\, ,\>=\frac{1}{p}\<\, ,\>_0$. The \dieu module $M(\ker \lambda)$
of the subgroup $\ker \lambda$ is equal to $M/M^t$. Hence
the condition $\ker \lambda\simeq \alpha_p\times \alpha_p$ is
equivalent to that $F$ and $V$ vanish on $M(\ker \lambda)=M/M^t$.

Since $\<\, ,\>$ is a perfect pairing on $FM_0$, that is,
$(FM_0)^t=FM_0$, 
we have \[ pM_0\subset M^t\subset FM_0\subset M\subset M_0. \]
Changing the notation, put $\ol M_0:=M_0/pM_0$ and let
\[ \<\, ,\>_0:\ol M_0\times \ol M_0\to k. \]
be the induced perfect pairing.
In $\ol M_0$, the subspace $\ol {M^t}$ is equal to  ${\ol M}^{\bot}$. 
Indeed,
\begin{equation}
  \begin{split}
    M^t & =\{m\in M_0; \<m,x\>_0\in pW\ \forall\, x\in M\},
    \\ \ol{M^t}& =\{m\in \ol M_0; \<m,x\>_0=0\ \forall\, x\in \ol
M\}={\ol M}^{\bot}.
  \end{split}
\end{equation}
From this we see that the condition $\ker\lambda\simeq \alpha_p\times
\alpha_p$ is equivalent to $\<\ol M, F\ol M\>=\<\ol M, V\ol M\>=0$.
Since $\ol {F M_0}=k<f_3,f_4>$, one has 
$\ol M=k <f_1', f_3,f_4>$ where $f_1'=af_1+bf_2$.
The condition $\<\ol M, F\ol M\>=\<\ol M, V\ol M\>=0$, same as
$\<f'_1,Ff'_1\>=\<f_1',Vf_1'\>=0$, gives the equation
$a^{p+1}+b^{p+1}=0$. This completes the proof. \qed
\end{proof}

Conversely, fix a polarized superspecial abelian surface
$(A,\lambda)$ such that $\ker \lambda \simeq \alpha_p\times \alpha_p$.
Then there are $p^2+1$ degree-$p$ isogenies $\varphi':(A_0,\lambda_0)\to
(A,\lambda)$ such that $A_0$ is superspecial and
$\varphi'^*\lambda=p\,\lambda_0$. Indeed, each isogeny $\varphi'$ always
has the property $\varphi'^*\lambda=p\,\lambda_0$ for a principal
polarization $\lambda_0$, and there are $|\bfP^1(\F_{p^2})|$ isogenies
with $A_0$ superspecial.

\subsection{}\label{s:45} We denote by $\calA'_{2,1,N,(p)}$ the moduli
space which parametrizes equivalence classes of 
isogenies $(\varphi':\ul A_0\to \ul A_1)$ of degree $p$, 
where $\ul A_1$ is an object in $\calA_{2,p,N}$ and $\ul A_0$ is an
object in $\calA_{2,1,N}$, such that
$\varphi'^* \lambda_1=p\,\lambda_0$ and $\varphi'_*\eta_0=\eta_1$.

There is a natural isomorphism from $\calA_{2,1,N,(p)}$
to $\calA'_{2,1,N,(p)}$. Given an object $(A,\lambda,\eta,H)$ in
$\calA_{2,1,N,(p)}$, let $\ul A_0:=\ul A$, $A_1:=A/H$ and
$\varphi':A_0\to A_1$ be the natural
projection. The polarization $p\,\lambda_0$ descends to one, 
denoted by $\lambda_1$, on $A_1$. 
Put $\eta_1:=\varphi'_*\eta_0$ and $\ul A_1=(A_1,\lambda_1,\eta_1)$.
Then $(\varphi':\ul A_0 \to \ul A_1)$ lies
in $\calA'_{2,1,N,(p)}$ and the morphism
\[ q: \calA_{2,1,N,(p)}\to \calA'_{2,1,N,(p)}, \quad
(A,\lambda,\eta,H)\mapsto (\varphi':\ul A_0 \to \ul A_1) \]
is an isomorphism.

We denote by $\calS'_{2,1,N,(p)}$ the supersingular locus of
$\calA'_{2,1,N,(p)}\otimes\Fpbar$. Thus we have
$\calS'_{2,1,N,(p)}\simeq \calS_{2,1,N,(p)}$. It is clear that
$S'_\xi\subset \calS'_{2,1,N,(p)}$ for each $\xi\in \Lambda$, 
and $S'_\xi\cap S'_{\xi'}=\emptyset$ if $\xi\neq \xi'$. 
For each $\gamma\in \Lambda_{2,1,N}$, 
let $S''_\gamma$ be the subspace of
$\calS'_{2,1,N,(p)}$ that consists of objects 
$(\varphi':\ul A_0 \to \ul A_1)$ with $\ul A_0=\ul A_\gamma$. One also 
has $S''_\gamma\cap S''_{\gamma'}=\emptyset$ if $\gamma\neq \gamma'$.

\begin{lm} \

{\rm (1)} Let $(M_0,\<\,,\>_0)$ be a separably quasi-polarized
  supersingular \dieu module of rank $4$ and suppose $a(M_0)=1$. Let
  $M_1:=(F,V)M_0$ and $N$ be the unique \dieu module containing $M_0$
  with $N/M_0=k$. Let $\<\,, \>_1:=\frac{1}{p}\<\,,\>_0$ be 
  the quasi-polarization for $M_1$.  
  Then one has $a(N)=a(M_1)=2$, $VN=M_1$, and 
  $M_1/M_1^t\simeq k\oplus k$ as \dieu modules.

{\rm (2)} Let $(M_1,\<\,,\>_1)$ be a quasi-polarized
  supersingular \dieu module of rank $4$. Suppose that $M_1/M_1^t$ is
  of length $2$, that is, the quasi-polarization has degree $p^2$. 
  \begin{itemize}
  \item [(i)] If $a(M_1)=1$, then letting $M_2:=(F,V)M_1$, one has that 
  $a(M_2)=2$ and $\<\,  ,\>_1$ is a separable quasi-polarization on $M_2$.  
  \item [(ii)] Suppose $(M_1,\<\,,\>_1)$ decomposes as the product of two 
  quasi-polarized \dieu submodules of rank $2$.
  Then there are a unique
  \dieu submodule $M_2$ of $M_1$ with $M_1/M_2=k$ and a unique \dieu
  module $M_0$ containing $M_1$ with $M_0/M_1=k$ so that
  $\<\,,\>_1$ {\rm (}resp.  $p\<\,,\>_1${\rm )} is a separable
  quasi-polarization on 
  $M_2$ {\rm(}resp. $M_0${\rm)}. 
  \item[(iii)] Suppose $M_1/M_1^t\simeq k\oplus k$ as \dieu modules. Let
  $M_2\subset M_1$ be any \dieu submodule with $M_1/M_2=k$, and
  $M_0\supset M_1$ be any \dieu overmodule with $M_0/M_1=k$. Then
  $\<\,,\>_1$ {\rm (}resp.  $p\<\,,\>_1${\rm )} is a separable 
  quasi-polarization on $M_2$ {\rm(}resp. $M_0${\rm)}.
  \end{itemize}
\end{lm}
This is well-known; the proof is elementary and omitted.

\begin{prop} \label{45} Notation as above.

{\rm (1)} One has
\[ \calS'_{2,1,N,(p)}=\left(\coprod_{\xi\in
    \Lambda} S'_\xi \right) \cup
\left(\coprod_{\gamma\in
    \Lambda_{2,1,N}}S''_\gamma\right). \]

{\rm (2)} The scheme $\calS'_{2,1,N,(p)}$ has ordinary double singular
points and 
\[ (\calS'_{2,1,N,(p)})^{\rm sing}=\left(\coprod_{\xi\in
    \Lambda} S'_\xi \right)\cap \left(\coprod_{\gamma\in
    \Lambda_{2,1,N}}S''_\gamma\right). \]
Moreover, one has
\[  |(\calS'_{2,1,N,(p)})^{\rm sing}|=
    |\Lambda_{2,1,N}|(p+1)=|\Lambda|(p^2+1). \]
\end{prop}
\begin{proof}
  (1) Let $(\varphi':\ul A_0 \to \ul A_1)$ be a point of
      $\calS'_{2,1,N,(p)}$. If $a(A_0)=1$, then $\ker \varphi'$ is the
      unique $\alpha$-subgroup of $A_0[p]$ and thus $\ul A_1\in
      \Lambda$. Hence this point lies in $S'_\xi$ for some $\xi$.
      Suppose that $\ul
      A_1$ is not in $\Lambda$, then there is a unique lifting
      $(\varphi_1':\ul A_0'\to \ul A_1)$ in $\calS'_{2,1,N,(p)}$ and
      the source $\ul A_0'$ is superspecial. Hence $\ul A_0=\ul A_0'$
      is superspecial and the
      point $(\varphi':\ul A_0 \to \ul A_1)$ lies in $S''_\gamma$ for
      some $\gamma$.

  (2) It is clear that the singularities only occur at the intersection
      of $S'_\xi$'s and $S''_\gamma$'s, as  $S'_\xi$ and $S''_\gamma$
      are smooth. Let $x=(\varphi':\ul A_\gamma
      \to \ul A_\xi)\in S'_\xi\cap S''_\gamma$. We know that the
      projection ${\rm pr}_0:\calS'_{2,1,N,(p)}\to \calS_{2,1,N}$ induces
      an isomorphism from $S'_\xi$ to $V'_\xi$. Therefore, ${\rm pr}_0$
      maps the one-dimensional subspace $T_x(S'_\xi)$ of
      $T_x(\calA'_{2,1,N,(p)}\otimes \Fpbar)$ onto the one-dimensional
      subspace $T_{{\rm pr}_0(x)}(V'_\xi)$ of $T_{{\rm pr}_0(x)}
      (\calA_{2,1,N}\otimes \Fpbar)$, where $T_x(X)$ denotes the tangent
      space of a variety $X$ at a point $x$. On the other hand, ${\rm
      pr}_0$ maps the subspace $T_x(S''_\gamma)$ to zero. This shows
      $T_x(S'_\xi)\neq T_x(S''_\gamma)$ in
      $T_x(\calA'_{2,1,N,(p)}\otimes \Fpbar)$; particularly
      $\calS'_{2,1,N,(p)}$ has ordinary double singularity at $x$.

      Since every singular point lies in both $S'_\xi$ and
      $S''_\gamma$ for some $\xi, \gamma$, by Subsection~\ref{s:44}
      each  $S'_\xi$ has $p^2+1$ singular points and each $S''_\gamma$
      has $p+1$ singular points. We get
\[  |(\calS'_{2,1,N,(p)})^{\rm sing}|=
    |\Lambda_{2,1,N}|(p+1),\quad\text{and}\quad
     |(\calS'_{2,1,N,(p)})^{\rm sing}|=|\Lambda|(p^2+1). \]
This completes the proof. \qed
\end{proof}
\begin{cor} \label{44} We have
\[ |(\calS'_{2,1,N,(p)})^{\rm sing}|=|\Sp_4(\Z/N\Z)|\cdot
\frac{(-1)\zeta(-1)\zeta(-3)}{4} (p-1)(p^2+1)(p+1) \]
and
\[ |\Lambda|=|\Sp_4(\Z/N\Z)|\cdot
\frac{(-1)\zeta(-1)\zeta(-3)}{4} (p^2-1). \]
\end{cor}
\begin{proof}
  This follows from Corollary~\ref{33} and (2) of
  Proposition~\ref{45}. \qed
\end{proof}

Note that the evaluation of $|\Lambda|$ here is different from that
given in Katsura-Oort \cite{katsura-oort:surface}. Their method does
not rely on the mass formula but the computation is more complicated.

Since $\calS_{2,1,N,(p)}\simeq\calS'_{2,1,N,(p)}$
    (Subsection~\ref{s:45}), Theorem~\ref{12} 
    follows from Proposition~\ref{45} and Corollary~\ref{44}.

As a byproduct, we obtain the description of the supersingular
locus $\calS_{2,p,N}$.

\begin{thm} \ 

 {\rm (1)} The scheme $\calS_{2,p,N}$ is equi-dimensional and
  each irreducible component is isomorphic to
  $\bfP^1$.

{\rm (2)} The scheme $\calS_{2,p,N}$  has $|\Lambda_{2,1,N}|$
irreducible components. 

{\rm (3)} The singular locus of $\calS_{2,p,N}$ consists of superspecial
points $(A,\lambda,\eta)$ with $\ker \lambda\simeq \alpha_p\times
\alpha_p$, and thus $|\calS_{2,p,N}^{\rm sing}|=|\Lambda|$. Moreover,
at each singular point there are $p^2+1$ irreducible
components passing through and intersecting transversely.

{\rm (4)} The natural morphism ${\rm pr}_1: \calS_{2,1,N,(p)}\to
 \calS_{2,p,N}$ contracts
 $|\Lambda|$ projective lines onto the singular locus of
$\calS_{2,p,N}$. 
\end{thm}

\section{The class numbers $H_n(p,1)$ and $H_n(1,p)$}
\label{sec:05}

In this section we describe the arithmetic part of the results in Li
and Oort \cite{li-oort}. Our references are Ibukiyama-Katsura-Oort
\cite[Section 2]{ibukiyama-katsura-oort} and Li-Oort \cite[Section
4]{li-oort}.

Let $B$ be the definite quaternion algebra over $\Q$ with discriminant
$p$, and $\calO$ be a maximal order of $B$.
Let $V=B^{\oplus n}$,  regarded as a left $B$-module of row vectors, and
let $\psi(x,y)=\sum_{i=1}^n x_i \bar y_i$ be the standard hermitian
form on $V$, where $y_i\mapsto \bar y_i$ is the canonical involution
on $B$. Let $G$ be the group of $\psi$-similitudes over $\Q$; its
group of $\Q$-points is
\[ G(\Q):=\{h\in M_n(B)\, |\, h \bar h^t=r I_n \text{ for some $r\in
  \Q^\times$}\, \}. \]

Two $\calO$-lattices $L$ and $L'$ in $B^{\oplus n}$ are called {\it
  globally
  equivalent} (denoted by $L\sim L'$) if $L'=L  h$ for some $h\in G(\Q)$.
For a finite place $v$ of $\Q$, we write $B_v:=B\otimes \Q_v$,
$\calO_v:=\calO\otimes \Z_v$ and $L_v:=L\otimes \Z_v$.
Two $\calO$-lattices $L$ and $L'$ in $B^{\oplus n}$ are called {\it
locally equivalent at $v$} (denoted by $L_v\sim L_v'$) if $L_v'=L_v  h_v$ for
  some $h_v\in G(\Q_v)$.
A {\it genus} of $\calO$-lattices is a set of (global) $\calO$-lattices in
$B^{\oplus n}$ which are equivalent to each other locally at every
finite place $v$.

Let
  \[ N_p=\calO_p^{\oplus n}\cdot
  \begin{pmatrix}
    I_{r} & 0 \\ 0 & \pi I_{n-r}
  \end{pmatrix}\cdot \xi\subset B_p^{\oplus n}, \]
where $r$ is the integer $[n/2]$, $\pi$ is a uniformizer in
$\calO_p$, and $\xi$ is an element in $\GL_n(B_p)$ such that
\[ \xi \bar \xi^t=\text{anti-diag$(1,1,\dots,1)$}. \]

\begin{defn}
(1) Let $\calL_n(p,1)$ denote the set of global equivalence classes of
    $\calO$-lattices $L$ in $B^{\oplus n}$ such that $L_v\sim
    \calO_v^{\oplus n}$ at every finite place $v$. The genus
    $\calL_n(p,1)$ is called the {\it principal genus}, and let
    $H_n(p,1):=|\calL_n(p,1)|$.

(2) Let $\calL_n(1,p)$ denote the set of global equivalence classes of
    $\calO$-lattices $L$ in $B^{\oplus n}$ such that $L_p\sim N_p$ and
    $L_v\sim \calO_v^{\oplus n}$ at every finite place $v\neq p$.
    The genus $\calL_n(p,1)$ is called the {\it non-principal genus},
    and let
    $H_n(1,p):=|\calL_n(1,p)|$.
\end{defn}

Recall (Section 3) that $\Lambda_g$ is the set of isomorphism classes of
$g$-dimensional principally polarized superspecial abelian varieties
over $\Fpbar$. When $g=2D>0$ is even, we denote by $\Lambda^*_{g,p^D}$ 
the set of isomorphism classes of $g$-dimensional polarized
superspecial abelian varieties $(A,\lambda)$ of degree $p^{2D}$
over $\Fpbar$ 
satisfying $\ker \lambda=A[F]$.

Let $\calA_{g,1}$ be the {\it coarse moduli scheme} of
$g$-dimensional principally polarized abelian varieties, and let
$\calS_{g,1}$ be the supersingular locus of $\calA_{g,1}\otimes \Fpbar$.
Recall (Subsection~\ref{sec:21}) that $\Pi_0(\calS_{g,1})$ denotes the
set of irreducible components of $\calS_{g,1}$.

\begin{thm}[Li-Oort]\label{52} We have
  \[ |\Pi_0(\calS_{g,1})|=
  \begin{cases}
    |\Lambda_g| &\text{if $g$ is odd;} \\
    |\Lambda^*_{g,p^D}| &\text{if $g=2D$ is even.}\\
  \end{cases} \]
\end{thm}

The arithmetic part for $\Pi_0(\calS_{g,1})$ is given by the following

\begin{prop}\label{53} \

  {\rm (1)} For any positive integer $g$, one has $|\Lambda_g|=H_g(p,1)$.

  {\rm (2)} For any even positive integer $g=2D$,
      one has $|\Lambda^*_{g,p^D}|=H_g(1,p)$.
\end{prop}
\begin{proof}
  (1) See \cite[Theorem 2.10]{ibukiyama-katsura-oort}. (2) See
      \cite[Proposition 4.7]{li-oort}.
\end{proof}

\def\Gr{\mathrm{Gr}}
\section{Correspondence computation} \label{sec:06}

\subsection{}
\label{sec:61}
Let $M_0$ be a superspecial \dieu module over $k$ of rank $2g$, and call
\[ \tilde M_0:=\{x\in M_0; F^2x=px\}, \] 
the skeleton of $M_0$ (cf. \cite[5.7]{li-oort}). We know that
 $\tilde M_0$ is a \dieu module over $\F_{p^2}$ and
 $\tilde M_0\otimes_{W(\F_{p^2})} W(k)=M_0$. The vector space 
$\tilde M_0/V \tilde M_0$ defines an $\F_{p^2}$-structure of the
 $k$-vector space $M_0/VM_0$.

Let $\Gr(n,m)$ be the Grassmannian variety of $n$-dimensional subspaces in
an $m$-dimensional vector space. Suppose $M_1$ is a \dieu submodule of
$M_0$ such that
\[ VM_0\subset M_1\subset M_0, \quad \dim_k M_1/VM_0=r, \]
for some integer $0\le r\le g$.
As $\dim M_0/VM_0=g$, the subspace $\ol M_1:=M_1/VM_0$ corresponds to
a point in $\Gr(r,g)(k)$.

\begin{lm}\label{61}
  Notation as above. Then $M_1$ is superspecial (i.e.
  $F^2M_1=pM_1$) if and only if $\ol M_1\in \Gr(r,g)(\F_{p^2})$.
\end{lm}
\begin{proof}
  If $M_1$ is generated by $V \tilde M_0$ and
  $x_1,x_2,\dots,x_r$, $x_i\in \tilde M_0$ over $W$. Then
  $\tilde M_1$ generates $M_1$ and thus
  $F^2M_1=pM_1$. Therefore, $M_1$ is superspecial.

Conversely if $M_1$ is superspecial, then we have
\[ V\tilde M_0\ \subset\  \tilde M_1
\ \subset \tilde M_0. \]
Therefore, $\tilde M_1$ gives rise to an element in
$\Gr(r,g)(\F_{p^2})$. \qed
\end{proof}

\subsection{}
Let $L(n,2n)\subset {\rm Gr}(n,2n)$ be the Lagrangian variety of
maximal isotropic subspaces 
in a $2n$-dimensional vector space with a non-degenerate
alternating form.

From now on $g=2D$ is an even positive integer. 
Recall (in Introduction and Section~\ref{sec:05}) that
$\Lambda^*_{g,p^{D},N}$ denotes the set of superspecial
(geometric) points $(A,\lambda,\eta)$ in
$\calA_{g,p^D,N}\otimes \Fpbar$ satisfying $\ker \lambda=A[F]$.

\begin{lm}\label{62} Let $(A_2,\lambda_2,\eta_2)\in
  \Lambda^*_{g,p^D,N}$ and $(M_2,\<\,,\>_2)$ 
be the associated \dieu module.
  There is a $W$-basis $e_1,\dots,
  e_{2g}$ for $M_2$ such that for $1\le i\le g$
  \[ Fe_i=e_{g+i},\quad Fe_{g+i}=pe_i, \]
  and the non-zero pairings are
  \[  \<e_i,e_{D+i}\>_2=-\<e_{D+i},e_i\>_2=\frac{1}{p}, \]     
\[  \<e_{g+i},e_{g+D+i}\>_2=-\<e_{g+D+i},e_{g+i}\>_2=1,\] 
for $1\le i\le D$.  
\end{lm}
\begin{proof}
  Use the same argument of Lemma~\ref{42} (2). 
\end{proof}

\subsection{}
Let $\calH$ be the moduli space over $\Z_{(p)}[\zeta_N]$ which 
parametrizes equivalence classes of objects $(\varphi:\ul A_1\to
\ul A_2)$,  where 
\begin{itemize}
\item $\ul A_1=(A_1,\lambda_1,\eta_1)$ is an object in
$\calA_{g,1,N}$,
\item $\ul A_2=(A_2,\lambda_2,\eta_2)$ is an object in
$\calA_{g,p^D,N}$, and
\item  $\varphi:A_1\to A_2$ is an
isogeny of degree $p^D$
satisfying $\varphi^*\lambda_2=p\lambda_1$ and
$\varphi_*\eta_1=\eta_2$.
\end{itemize}
 
The moduli space $\calH$ with two natural projections gives the
following correspondence:
\[ \xymatrix{
 &  \calH \ar[ld]_{\pr_1} \ar[rd]^{\pr_2} & \\
\calA_{g,1,N} & & \calA_{g,p^D,N}.
} \]

Let $\calS$ be the supersingular locus of $\calH\otimes \Fpbar$, which 
is the reduced closed subscheme consisting of supersingular points
(either $A_1$ or $A_2$ is supersingular, or equivalently both are
supersingular). 
Restricting the natural projections on $\calS$, we have the following
correspondence  
\[ \xymatrix{
 &  \calS \ar[ld]_{\pr_1} \ar[rd]^{\pr_2} & \\
\calS_{g,1,N} & &\calS_{g,p^D,N}.
} \]

Suppose that $\ul A_2 \in\Lambda^*_{g,p^D,N}$. Let $(\varphi:\ul
A_1\to \ul A_2)\in \calS(k)$ be a point in the pre-image
$\pr^{-1}_2(\ul A_2)$, and let
$(M_1,\<\,,\>_1)$ be the \dieu module associated to $\ul A_1$. We
have
\[ M_1\subset M_2,\quad p \<\,,\>_2=\<\,,\>_1, \quad FM_2=M_2^t. \]
Since $\ul A_2$ is superspecial and $ \<\,,\>_1$ is a perfect pairing 
on $M_1$, we get
\[ FM_2=VM_2=M_2^t, \quad M^t_2\subset M_1^t=M_1. \]
Therefore, we have
\[ FM_2=VM_2\subset M_1\subset M_2, \quad \dim_k M_1/VM_2=D. \]
Put $\<\,,\>:=p\<\,,\>_2$. The pairing
\[ \<\,,\>: M_2\times M_2\to W \]
induces a pairing
\[ \<\,,\>: M_2/VM_2\times M_2/VM_2\to k \]
which is perfect (by Lemma~\ref{62}). Furthermore, $M_1/VM_2$ is a
maximal isotropic 
subspace for the pairing $ \<\,,\>$. This is because $\<\,,\>_1$ is a
perfect pairing on $M_1$ and $\dim M_1/VM_2=D$ is the maximal
dimension of isotropic subspaces. We conclude that the point
$(\varphi:\ul A_1\to \ul A_2)$ lies in $\pr^{-1}_2(\ul A_2)$ if and
only if
$VM_2\subset M_1\subset M_2$ and $M_1/VM_2$ is a maximal isotropic 
subspace of the symplectic space $(M_2/VM_2,\<\, ,\>)$. By
Lemma~\ref{61}, we have proved

\begin{prop}\label{63} Let $\ul A_2$ be a point in $\Lambda^*_{g,p^D,N}$.

  {\rm (1)} The pre-image $\pr^{-1}_2(\ul A_2)$
      is naturally isomorphic to the projective variety $L(D,
      2D)$ over $k$.

  {\rm (2)} The set $\pr_1^{-1} (\Lambda_{g,1,N})\cap \pr^{-1}_2(\ul
      A_2)$ is in bijection with $L(D, 2D)(\F_{p^2})$, where the
      $W(\F_{p^2})$-structure of $M_2$ is given by
      the skeleton $\tilde M_2$.
\end{prop}

\subsection{} \label{sec:64}
We compute $\pr_1^{-1}(\ul A_1) \cap \pr_2^{-1} 
(\Lambda^*_{g,p^D,N})$ for a point 
$\ul A_1$ in $\Lambda_{g,1,N}$. 
Let $\calT$ be the closed subscheme of $\calS$ consisting of
the points $(\varphi:\ul A_1\to \ul A_2)$ such that $\ker
\lambda_2=A_2[F]$. We compute the closed subvariety 
$\pr_1^{-1}(\ul A_1) \cap \calT$ first. 

Let $\ul A_1\in \Lambda_{g,1,N}$, and let $(\varphi:\ul
A_1\to \ul A_2)\in \calS(k)$ be a point in the pre-image
$\pr^{-1}_1(\ul A_1)$.
Let $(M_1,\<\,,\>_1)$ and $(M_2,\<\,,\>_2)$ be the
\dieu modules associated to $\ul A_1$ and $\ul A_2$, respectively. 

One has
\[ M_1^t=M_1\supset M_2^t=FM_2, \]
and thus has
\[ FM_2\subset M_1\subset M_2\subset p^{-1}VM_1. \]
Since $M_1$ is superspecial, $p^{-1}VM_1=p^{-1}FM_1$. Put
$M_0:=p^{-1}VM_1$ and $\<\,,\>:=p\<\,,\>_2$. We have
\[ pM_0\subset M_2^t=FM_2\subset M_1=VM_0\subset M_2\subset M_0 \]
and that $\<\,,\>$ is a perfect pairing on $M_0$. By Proposition 6.1 of
\cite{li-oort} (cf. Lemma~\ref{42} (1)), there is a $W$-basis
$f_1\dots,f_{2g}$ for $M_0$ such that for $1\le i\le g$
\[ Ff_{i}=f_{g+i},\quad Ff_{g+i}=pf_i, \]
and the non-zero pairings are
\[ \<f_i, f_{g+i}\>=-\<f_{g+i}, f_i\>=\beta_1, \quad \forall\, 1\le i\le g\]
where $\beta_1\in W(\F_{p^2})^\times$ with $\beta_1^\sigma=-\beta_1. $
In the vector space $\ol M_0:=M_0/pM_0$, $\ol M_2$ is a vector
subspace over $k$ of dimension $g+D$ with
\[ \ol M_2\supset \ol {VM}_0=k<f_{g+1},\dots,f_{2g}>
\quad\text{and\ \  } \<\ol M_2,
F\ol M_2\>=0. \]
We can write
\[ \ol M_2=k<v_1,\dots,v_D>+\ol {VM}_0, \quad v_i=\sum_{r=1}^g
a_{ir} f_r. \]
One computes
\[ \<v_i, Fv_j\>=\<\sum_{r=1}^g
a_{ir} f_r, F(\sum_{q=1}^g a_{jq} f_q)\>= \<\sum_{r=1}^g
a_{ir} f_r, \sum_{q=1}^g a_{jq}^p f_{g+q}\>=\beta_1 \sum_{r=1}^g
a_{ir}  a_{jr}^p. \]
This computation leads us to the following definition.

\subsection{}
Let $V:=\F_{p^2}^{2n}$. For any field $K\supset \F_{p^2}$, we put
$V_K:=V\otimes_{\F_{p^2}} K$ and define a pairing on $V_K$
\[ \<\,,\>':V_K\times V_K\to K,\quad \<(a_i),(b_i)\>':=\sum_{i=1}^{2n}\,
a_i\, b_i^p. \]
Let $\bfX(n,2n) \subset {\rm Gr}(n,2n)$ be the subvariety over $\F_{p^2}$
which parametrizes $n$-dimensional (maximal) isotropic subspaces in $V$
with respect to the pairing $\<\,,\>'$. 

With the computation in Subsection~\ref{sec:64} and Lemma~\ref{61}, we
have proved
\begin{prop}\label{64}
  Let $\ul A_1$ be a point in $\Lambda_{g,1,N}$ and $g=2D$. 

  {\rm (1)} The intersection $\pr_1^{-1}(\ul A_1)\cap \calT$ is
      naturally isomorphic to the projective variety $\bfX(n,2n)$ over
      $k$.

  {\rm (2)} The set 
      $\pr_1^{-1}(\ul A_1)\cap \pr_2^{-1}(\Lambda^*_{g,p^D,N})$ is 
      in bijection with $\bfX(D,2D)(\F_{p^2})$.
\end{prop}

When $g=2$, Proposition~\ref{64} (2) is a result of Moret-Bailly
(Lemma~\ref{43}).

\subsection{}
\label{sec:66}

To compute $\bfX(n,2n)(\F_{p^2})$, we show that it is 
the set of rational points of a homogeneous space under the
quasi-split group $U(n,n)$.

Let $V=\F_{p^2}^{2n}$ and let $x\mapsto \bar x$ be the involution of
$\F_{p^2}$ 
over $\Fp$. Let $\psi((x_i),(y_i))=\sum_i x_i\bar y_i$ be the standard
hermitian form on $V$. For any field $K\supset\Fp$, we put
$V_K:=V\otimes_{\Fp}K$ and extend $\psi$ to a from
\[ \psi:V_K\otimes V_K\to \F_{p^2}\otimes K \]
by $K$-linearity.

Let $U(n,n)$ be the group of automorphisms of $V$ that preserve the
hermitian form $\psi$.
Let $LU(n,2n)(K)$ be the space of $n$-dimensional (maximal) isotropic
$K$-subspaces in $V_K$ with respect to $\psi$. We know that $LU(n,2n)$ is
a projective scheme over $\F_p$ of finite type, and this is a
homogeneous space under 
$U(n,n)$. It follows from the definition that
\[ LU(n,2n)(\Fp)=\bfX(n,2n)(\F_{p^2}). \]
However, the space $LU(n,2n)$ is not isomorphic to
the space $\bfX(n,2n)$ over $k$.

Let $\wt \Lambda$ be the subset of $\calS$ consisting of elements
$(\varphi:\ul A_1\to \ul A_2)$ such that $\ul A_2\in
\Lambda^*_{g,p^D,N}$ and $\ul A_1\in \Lambda_{g,1,N}$. We have
natural projections
\[ \xymatrix{
 &  \wt \Lambda \ar[ld]_{\pr_1} \ar[rd]^{\pr_2} & \\
\Lambda_{g,1,N} & & \Lambda^*_{g,p^D,N}.
} \]

By Propositions~\ref{63} and \ref{64}, and Subsection~\ref{sec:66}, we
have proved

\begin{prop}\label{65}
Notation as above, one has
\[ |\wt \Lambda|=|L(D,2D)(\F_{p^2})|\cdot
|\Lambda^*_{g,p^D,N}|=|LU(D,2D)(\Fp)|\cdot
 |\Lambda_{g,1,N}|. \]
\end{prop}

\begin{thm}\label{66}
  We have
\[ |\Lambda^*_{g,p^D,N}|=
  |\Sp_{2g}(\Z/N\Z)|\cdot
  \frac{(-1)^{g(g+1)/2}}{2^g} \left \{ \prod_{i=1}^g \zeta(1-2i)
  \right \}\cdot \prod_{i=1}^{D}(p^{4i-2}-1). \]
\end{thm}
\begin{proof}
  We compute in Section~\ref{sec:07} that
\[ |L(D,2D)(\F_{p^2})|= \prod_{i=1}^D (p^{2i}+1), \]
\[ |LU(D,2D)(\Fp)|=\prod_{i=1}^D (p^{2i-1}+1). \]
Using Proposition~\ref{65}, Theorem~\ref{31} and Lemma~\ref{32}, we
get the value of $|\Lambda^*_{g,p^D,N}|$. \qed
\end{proof}

\subsection{Proof of Theorem~\ref{13}}
\label{sec:67} 
By a theorem of Li and Oort (Theorem~\ref{52}), we know
\[ |\Pi_0(\calS_{g,1,N})|=
\begin{cases}
 |\Lambda_{g,1,N}| &  \text{if $g$ is odd;} \\
 |\Lambda^*_{g,p^{D},N}| &  \text{if $g=2D$ is even.}
\end{cases} \]
Note that the result of Li and Oort is formulated for the coarse
moduli space $\calS_{g,1}$. However, it is clear that adding
the level-$N$ structure yields a modification as above. 
Theorem~\ref{13} then follows from Theorem~\ref{31}, Lemma~\ref{32} and
Theorem~\ref{66}. \qed

\section{$L(n,2n)(\F_q)$ and $LU(n,2n)(\F_q)$}
\label{sec:07}
Let $L(n,2n)$ be the Lagrangian variety of maximal isotropic subspaces
in a $2n$-dimensional vector space $V_0$ with a non-degenerate
alternating form $\psi_0$.

\begin{lm}\label{71}
$ |L(n,2n)(\Fq)|= \prod_{i=1}^n (q^i+1).$
\end{lm}
\begin{proof}
Let $e_1,\dots,e_{2n}$ be the standard symplectic basis for $V_0$. The
group $\Sp_{2n}(\Fq)$ acts transitively on the space $L(n,2n)(\Fq)$.
For $h\in \Sp_{2n}(\Fq)$, the map $h\mapsto \{he_1,\dots he_{2n}\}$
induces a bijection between $\Sp_{2n}(\Fq)$ and the set $\calB(n)$ of
ordered symplectic bases $\{v_1,\dots,v_{2n}\}$ for $V_0$.
The first vector $v_1$ has $q^{2n}-1$ choices. The first companion 
vector $v_{n+1}$ has $(q^{2n}-q^{2n-1})/(q-1)$ choices as it 
does not lie in the hyperplane $v_1^{\bot}$ and we require 
$\psi_0(v_1,v_{n+1})=1$.
The remaining ordered symplectic basis can be chosen from the
complement 
$\Fq<v_1,v_{n+1}>^{\bot}$. Therefore, we have proved the recursive formula
\[ |\Sp_{2n}(\Fq)|=(q^{2n}-1)q^{2n-1}|\Sp_{2n-2}(\Fq)|. \]
From this, we get
\[ |\Sp_{2n}(\Fq)|=q^{n^2}\prod_{i=1}^n (q^{2i}-1). \]
Let $P$ be the stabilizer of the standard maximal isotropic subspace
$\Fq <e_1,\dots,e_n>$. It is easy to see that
\[ P=\left\{
    \begin{pmatrix}
      A & B \\ 0 & D \\
    \end{pmatrix} ; AD^t=I_n,\ BA^t=AB^t \right\}. \]
The matrix $BA^t$ is symmetric and the space of $n\times n$
symmetric matrices has dimension $(n^2+n)/2$.
This yields
\[ |P|=q^{\frac{n^2+n}{2}}\, |\GL_n(\Fq)|=q^{n^2} 
\prod_{i=1}^n (q^i-1) \] as one has
\[ |\GL_n(\Fq)|=q^{\frac{n^2-n}{2}} \prod_{i=1}^n (q^i-1). \]
Since $L(n,2n)(\Fq)\simeq \Sp_{2n}(\Fq)/P$, we get
$ |L(n,2n)(\Fq)|= \prod_{i=1}^n (q^i+1).$ \qed
\end{proof}

\subsection{}
\label{sec:71}
Let $V=\F_{q^2}^{2n}$ and let $x\mapsto \bar x=x^q$ be the involution of
$\F_{q^2}$ over $\Fq$. Let $\psi((x_i),(y_i))=\sum_i x_i\bar y_i$
be the standard hermitian form on $V$. For any field $K\supset\Fq$,
we put
$V_K:=V\otimes_{\Fp}K$ and  extend $\psi$ to a
from
\[ \psi:V_K\otimes V_K\to \F_{q^2}\otimes_{\Fq} K \]
by $K$-linearity.

Let $U(n,n)$ be the group of automorphisms of $V$ that preserve the
hermitian form $\psi$.
Let $LU(n,2n)(K)$ be the set of $n$-dimensional (maximal) isotropic
$K$-subspaces in $V_K$ with respect to $\psi$. We know that $LU(n,2n)$ is
a homogeneous space under $U(n,n)$.
  Let \[ I_m:=\{\ul a=(a_1,\dots, a_m)\in \F_{q^2}^m\,;\
  Q(\ul a)=0 \}, \]
  where $Q(\ul a)=a_1^{q+1}+\dots+a_m^{q+1}$.
\begin{lm}\label{72} We have
   $|I_m|=q^{2m-1}+(-1)^m q^m+(-1)^{m-1}q^{m-1}$.
\end{lm}
\begin{proof}
For $m>1$, consider the projection $p:I_m\to \F_{q^2}^{m-1}$ which
sends $(a_1,\dots,a_m)$ to $(a_1,\dots,a_{m-1})$. Let $I^c_{m-1}$
be the complement of $I_{m-1}$ in $\F_{q^2}^{m-1}$. If $x\in
I_{m-1}$, then the pre-image $p^{-1}(x)$ consists of one element.
If $x\in I^c_{m-1}$, then the pre-image $p^{-1}(x)$ consists
of solutions of the equation $a_m^{q+1}=-Q(x)\in \Fq^\times$ and
thus $p^{-1}(x)$ has $q+1$ elements. Therefore,
$|I_m|=|I_{m-1}|+(q+1)|I_{m-1}^c|$. From this we get the recursive
formula 
\[ |I_m|=(q+1)q^{2(m-1)}-q|I_{m-1}|. \]
We show the lemma by induction. When $m=1$, $|I_m|=1$ and the
statement holds. Suppose the statement holds for $m=k$, i.e.
$|I_k|=q^{2k-1}+(-1)^k q^k+(-1)^{k-1}q^{k-1}$. When $m=k+1$,
\begin{equation*}
\begin{split}
|I_{k+1}|&=(q+1)q^{2k}-q[q^{2k-1}+(-1)^k q^k+(-1)^{k-1}q^{k-1}] \\
&=q^{2k+1}+(-1)^{k+1} q^{k+1}+(-1)^{k}q^{k}. \\
\end{split}
\end{equation*}
This completes the proof. \qed
\end{proof}

\begin{prop}\label{73}
  $|LU(n,2n)(\Fq)|=\prod_{i=1}^n (q^{2i-1}+1). $
\end{prop}
\begin{proof}
  We can choose a new basis $e_1,\dots,e_{2n}$ for $V$ such that 
  the non-zero pairings are
\[ \psi(e_i,e_{n+i})=\psi(e_{n+i},e_i)=1, \quad \forall\, 1\le i\le
  n. \]
The representing matrix for $\psi$ with respect to
$\{e_1,\dots,e_{2n}\}$ is
\[ J=
\begin{pmatrix}
  0 & I_n \\ I_n & 0
\end{pmatrix}. \]
Let $P$ be the stabilizer of the standard maximal isotropic subspace
$\F_{q^2} <e_1,\dots,e_n>$. It is easy to see that
\[ P=\left\{
    \begin{pmatrix}
      A & B \\ 0 & D \\
    \end{pmatrix} ; AD^*=I_n,\ BA^*+AB^*=0 \right\}. \]
The matrix $BA^*$ is skew-symmetric hermitian. The space of $n\times n$
skew-symmetric hermitian matrices has dimension $n^2$ over
$\Fq$. Indeed, the diagonal consists of 
entries in the kernel of the trace; this gives dimension $n$. 
The upper triangular has $(n^2-n)/2$ entries in
$\F_{q^2}$; this gives dimension $n^2-n$. Hence,
\[ |P|=q^{n^2}|\GL_n(\F_{q^2})|=q^{2n^2-n} \prod_{i=1}^n (q^{2i}-1).\]

We compute $|U(n,n)(\Fq)|$. For $h\in U(n,n)(\Fq)$, the map
\[ h\mapsto \{he_1,\dots, he_{2n}\} \]
gives a bijection between $U(n,n)(\Fq)$ and the set $\calB(n)$ of
ordered bases $\{v_1,\dots,v_{2n}\}$ for which the representing matrix
of $\psi$ is $J$. The first vector $v_1$ has
\[ |I_{2n}|-1=q^{4n-1}+q^{2n}-q^{2n-1}-1=(q^{2n}-1)(q^{2n-1}+1). \]
choices (Lemma~\ref{72}). For the choices of the companion vector
$v_{n+1}$ with $\psi(v_{n+1}, v_{n+1})=0$ and $\psi(v_1, v_{n+1})=1$,
consider the set
\[Y:=\{v\in V; \psi(v_1,v)=1\}. \]
Clearly, $|Y|=q^{4n-2}$. The additive group $\F_{q^2}$ acts on $Y$ by
$a\cdot v=v+av_1$ for $a\in \F_{q^2}$, $v\in Y$. It follows from
\[ \psi(v+a v_1,v+av_1)=\psi(v,v)+\bar a+a \]
that every orbit $O(v)$ contains an isotropic vector $v_0$ and any
isotropic vector in $O(v)$ has the form $v_0+av_1$ with
$\bar a+a=0$. Hence, the vector $v_{n+1}$ has
\[ \frac{|\,Y|\, q}{q^2}=q^{4n-3} \]
choices. In conclusion, we have proved the recursive formula
\[ |U(n,n)(\Fq)|=(q^{2n}-1)(q^{2n-1}+1)q^{4n-3} |U(n-1,n-1)(\Fq)|. \]
It follows that
\[ |U(n,n)(\Fq)|=q^{2n^2-n} \prod_{i=1}^n (q^{2i}-1)(q^{2i-1}+1). \]
Since $LU(n,2n)(\Fq)\simeq U(n,n)(\Fq)/P$, we get
$|LU(n,2n)(\Fq)|=\prod_{i=1}^n(q^{2i-1}+1)$.\qed

\end{proof}


\end{document}